\documentclass[11pt, a4paper]{amsart}
\usepackage{amsmath,amsfonts,amssymb,amsthm}
\usepackage{graphicx} 
\usepackage{enumitem}

\title{A note on the Erd\H{o}s Matching Conjecture}
\author{Ryan R. Martin}
\address{Iowa State University, Ames, Iowa, USA} 
\thanks{Research partially supported by Simons Foundation Collaboration Grant for Mathematicians \#709641}
\email{rymartin@iastate.edu}
\author{Bal\'azs Patk\'os}
\address{HUN-REN Alfr\'ed R\'enyi Institute of Mathematics} 
\thanks{Research supported by NKFIH under grant FK132060} 
\email{patkosb@gmail.com}
\date{March 2024}

\newtheorem{thm}{Theorem}
\newtheorem{cor}[thm]{Corollary}
\newtheorem{lem}[thm]{Lemma}
\newtheorem{prop}[thm]{Proposition}
\newtheorem{conjecture}[thm]{Conjecture}

\newcommand{\F}{{\mathcal F}}
\newcommand{\cF}{{\mathcal F}}
\newcommand{\cS}{{\mathcal S}}
\newcommand{\cG}{{\mathcal G}}
\newcommand{\cA}{{\mathcal A}}
\newcommand{\cB}{{\mathcal B}}
\newcommand{\T}{{\mathcal T}}

\begin{document}

\begin{abstract}
    The Erd\H os Matching Conjecture states that the maximum size $f(n,k,s)$ of a family $\cF\subseteq \binom{[n]}{k}$ that does not contain $s$ pairwise disjoint sets is $\max\{|\cA_{k,s}|,|\cB_{n,k,s}|\}$, where $\cA_{k,s}=\binom{[sk-1]}{k}$ and $\cB_{n,k,s}=\{B\in \binom{[n]}{k}:B\cap [s-1]\neq \emptyset\}$. The case $s=2$ is simply the Erd\H os-Ko-Rado theorem on intersecting families and is well understood. The case $n=sk$ was settled by Kleitman and the uniqueness of the extremal construction was obtained by Frankl. Most results in this area show that if $k,s$ are fixed and $n$ is large enough, then the conjecture holds true. Exceptions are due to Frankl who proved the conjecture and considered variants for $n\in [sk,sk+c_{s,k}]$ if $s$ is large enough compared to $k$. A recent manuscript by Guo and Lu considers non-trivial families with matching number at most $s$ in a similar range of parameters.

    In this short note, we are concerned with the case $s\ge 3$ fixed, $k$ tending to infinity and $n\in\{sk,sk+1\}$. For $n=sk$, we show the  stability of the unique extremal construction of size $\binom{sk-1}{k}=\frac{s-1}{s}\binom{sk}{k}$ with respect to minimal degree. As a consequence we derive  $\lim\limits_{k\rightarrow \infty}\frac{f(sk+1,k,s)}{\binom{sk+1}{k}}<\frac{s-1}{s}-\varepsilon_s$ for some positive constant $\varepsilon_s$ which depends only on $s$.
\end{abstract}
\maketitle

\section{Introduction}
We use standard notation. Let $[n]$ denote the set of the first $n$ positive integers. For a set $X$, $\binom{X}{k}$ stands for the family of all $k$-element subsets of $X$ and a \textit{$k$-graph} (or $k$-uniform hypergraph) on $X$ is a subset of $\binom{X}{k}$. The members of $X$ are called \textit{vertices} and the $k$-sets are called \textit{$k$-edges} or \textit{hyperedges} or just \textit{edges}. A $k$-graph is \textit{intersecting} if every pair of edges has a nonempty intersection. A $k$-graph on $X$ for which there is an $x\in X$ that is in every edge is called a \textit{star}.

One of the major unsolved problems in extremal set theory is due to Erd\H os: The \textit{matching number} $\nu(\cF)$ of a family $\cF$ of sets is the size of the largest set of pairwise disjoint sets that $\cF$ contains. Let $f(n,k,s)$ denote the maximum size of a $k$-graph on $n$ vertices with matching number strictly smaller than $s$. Because $\nu(\cF)=1$ if and only if $\cF$ is intersecting, the well-known Erd\H{o}s-Ko-Rado theorem \cite{EKR} gives 
\begin{align*}
    f(n,k,2) = \begin{cases}
                \binom{n-1}{k-1}, &\mbox{ if $n\geq 2k$;} \\
                \binom{n}{k}, &\mbox{ if $n\leq 2k-1$.}
    \end{cases}
\end{align*}

The special case of $f(sk,k,s)$ was proved by Kleitman.
\begin{thm}[Kleitman~\cite{K}]\label{K}
    If $k\geq 2$ and $s\geq 2$ are positive integers, then $f(sk,k,s)=\binom{sk}{k}-\binom{sk-1}{k-1}=\frac{s-1}{s}\binom{sk}{k}$ with equality if and only if there is an element $x\in[n]$ such that every $k$-edge fails to contain $x$.
\end{thm}
There are two natural families of $k$-sets with matching number $s-1$:

\[\cA_{k,s}=\binom{[sk-1]}{k} \hskip 1truecm
\cB_{n,k,s}=\left\{B\in \binom{[n]}{k}:B\cap [s-1]\neq \emptyset\right\}.
\]
Erd\H os conjectured that for any values of $n,k,s$, one of these families achieves the maximum value.

\begin{conjecture}[Erd\H os Matching Conjecture (EMC) \cite{E}]
   For all $n,k,s$ with $n\ge sk$, we have $f(n,k,s)=\max\{|\cA_{k,s}|,|\cB_{n,k,s}|\}$. \label{conj:EMC}
\end{conjecture}

    It can be computed that $|\cA_{k,s}|<|\cB_{n,k,s}|$ for all $n\geq s(k+1)$. Most results concerning Conjecture~\ref{conj:EMC} determine that $f(n,k,s)=|\cB_{n,k,s}|$ if $n$ is large enough compared to $k$ and $s$ \cite{FLM,HLS}, with the current best bounds due to Frankl~\cite{Fra1.5} for $n\ge (2s+1)k-s$ and Frankl and Kupavskii~\cite{FK} for $n\ge \frac{5}{3}sk-\frac{2}{3}s$ if $s$ is large enough. The case $k=3$ is also completely settled by Frankl~\cite{Fra2}. Minimum degree versions of the problem are also studied \cite{B,GLJ,HZ} and general inequalities on $f(n,k,s)$ are known \cite{FK,H} for $n=cks$ with $c>1$ and $k$ and/or $s$ is large enough.

    Not too many results are known in the region when $n$ is very close to $sk$: Frankl~\cite{Fra3} proved that $f(n,k,s)=|\cA_{k,s}|$, if $sk\le n\le s(k+\varepsilon)$ where $\varepsilon$ depends only on $k$, but $s$ has to be large enough with respect to $k$ ($s\geq k^2+k+1$ suffices). This result was extended by Kolupaev and Kupavskii \cite{KK}. Frankl~\cite{Fra5} and recently Guo and Lu~\cite{GL} considered non-trivial families (i.e. those with no isolated vertices) in a similar range of parameters.
    
    Our aim is to obtain bounds on $f(n,k,s)$ when $n=sk+1$, but with $s$ fixed and $k$ being large. Observe that the EMC states that there should be a huge difference between the cases $s=2$ and $s\ge 3$: on the one hand $\frac{f(2k,k,2)}{\binom{2k}{k}}=\frac{1}{2}$ and $\frac{f(2k+1,k,2)}{\binom{2k+1}{k}}\rightarrow \frac{1}{2}$ as $k$ tends to infinity. On the other hand, $\frac{f(sk,k,s)}{\binom{sk}{k}}=\frac{s-1}{s}$, while $\frac{|\cA_{k,s}|}{\binom{sk+1}{k}}\rightarrow (\frac{s-1}{s})^2$ and $\frac{|\cB_{sk+1,k,s}|}{\binom{sk+1}{k}}\rightarrow 1-(1-\frac{1}{s})^{s-1}$ as $k$ tends to infinity, thus the EMC would yield an immediate drop in the limiting constant. Our main result shows that there is indeed a gap between $\frac{f(sk,k,s)}{\binom{sk}{k}}=\frac{s-1}{s}$ and $\frac{f(sk+1,k,s)}{\binom{sk+1}{k}}$.

\begin{thm}\label{drop}
    For any $s\ge 3$ there exists a positive real $\varepsilon_s$ such that $f(sk+1,k,s)\le \bigl(\frac{s-1}{s}-\varepsilon_s\bigr)\binom{sk+1}{k}$ holds for all $k$.
\end{thm}

To prove Theorem \ref{drop}, as an intermediate step, we will show the following stability version of Kleitman's bound on $f(sk,k,s)$ with respect to the minimum degree.

\begin{thm}\label{stab}
    There exist absolute constants $C$ and $\delta_0$ such that, if $s\ge 3$ and $\delta\le\delta_0$, then any family $\cF\subseteq \binom{[sk]}{k}$ with $\nu(\cF)\le s-1$ and minimum degree at least $\delta\binom{sk-1}{k-1}$ satisfies $|\cF|\le \bigl(\frac{s-1}{s}-\frac{(s-2)\delta}{s^3(s-1)C}\bigr)\binom{sk}{k}$.
\end{thm}

\section{Tools of the proofs: removal lemmas and shifting}

In this section, we introduce some terminology and results from the literature that we will use in our proofs of Theorem \ref{drop} and Theorem \ref{stab}. \textit{Removal lemmas} are widely used in combinatorics, they state that if some combinatorial structure $\cS$ contains only a few copies of some pattern $P$, then $\cS$ can be made $P$-free while almost keeping the structure $\cS$. We will need a removal lemma for intersecting families. A general result is the following.

\begin{thm}[Friedgut, Regev~\cite{FR}]\label{fr}
    Let $\gamma>0$ and $n,k$ positive integers with $\gamma n<k<(1/2-\gamma)n$. Then for any $\varepsilon>0$ there exists $\delta>0$ such that any $\cF\subseteq \binom{[n]}{k}$ with at most $\delta |\cF|\binom{n-k}{k}$ disjoint pairs can be made intersecting by removing at most $\varepsilon\binom{n-1}{k-1}$ sets from $\cF$.
\end{thm}

As we will point out in the last section, Theorem \ref{fr} could also be used to derive our results, but the following statement is better suited for our purposes.

\begin{thm}[Das, Tran~\cite{DT}]
    There exists an absolute constant $C$ such that if $n,k,\ell$ satisfy $n>2k\ell^2$, then for any $\cF\subseteq \binom{[n]}{k}$ of size $(\ell-\alpha)\binom{n-1}{k-1}$ with at most $\bigl(\binom{\ell}{2}+\beta\bigr)\binom{n-1}{k-1}\binom{n-k-1}{k-1}$ disjoint pairs where $\max\{2\ell|\alpha|,|\beta|\}\le \frac{n-2k}{(20C)^2n}$ there exists a family $\cS$ that is a union of $\ell$ stars such that $|\cF\triangle \cS|\le C\bigl((2\ell-1)\alpha+2\beta\bigr)\frac{n}{n-2k}\binom{n-1}{k-1}$.
\end{thm}

We state the case where $\ell=1,\alpha=0, n=sk>2k$ as a corollary, since this is what we will use in our proofs.

\begin{cor}\label{dp}
    There exists an absolute constant $C$ such that for any $\cF\subseteq \binom{[sk]}{k}$ of size $\binom{sk-1}{k-1}$ with at most $\frac{\beta}{s-1}\binom{sk-1}{k-1}\binom{(s-1)k}{k}$ disjoint pairs where $\beta\le \frac{s-2}{s(20C)^2}$ there exists a star $\cS$ such that $|\cF\triangle \cS|\le \frac{2s}{s-2}C\beta\binom{sk-1}{k-1}$.

    Consequently, if the maximum degree $\Delta(\cF)$ of a family $\cF\subseteq \binom{[sk]}{k}$ of size $\binom{sk-1}{k-1}$ is at most $(1-\delta)\binom{sk-1}{k-1}$ where $\delta\le \frac{1}{200C}$, then $\cF$ contains at least $\frac{\delta(s-2)}{2Cs(s-1)}\binom{sk-1}{k-1}\binom{(s-1)k}{k}$ disjoint pairs.
\end{cor}

Let us now turn to our other tool in proving our results. The \textit{shifting operation} was introduced by Erd\H{o}s, Ko, and Rado~\cite{EKR} and is a very frequently used powerful tool in extremal set theory. For a set $F$ and family $\cF$ of sets we define

\begin{align*}
    S_{i,j}(F) = \begin{cases}
               F\setminus \{j\}\cup \{i\}, &\mbox{ if $i\notin F, j\in F, F\setminus \{j\}\cup \{i\}\notin F$;} \\
                F, &\mbox{otherwise.}
    \end{cases}
\end{align*}
$S_{i,j}(\cF)=\{S_{i,j}(F):F\in \cF\}$. A family $\cF\subseteq \binom{[n]}{k}$ is \textit{left-compressed} if $S_{i,j}(\cF)=\cF$ for all $1\le i<j\le n$.
We will use the following lemma~\cite[Lemma 4.2(iv)]{Fra}.

\begin{lem}[Frankl~\cite{Fra}]\label{shiftnu}
    If $\F'\subseteq \binom{[n]}{k}$ has matching number $\nu(\cF)\le s$, then there is a left-compressed $\cF\subseteq \binom{[n]}{k}$ with $\nu(\cF)\le s$ and $\left|\F\right|=\left|\F'\right|$.
\end{lem}

We will need some more properties of left-compressed families. To state them, we need some definitions: for a family $\cF\subseteq \binom{[n]}{k}$ and $x,y\in [n]$, we write $\cF_x=\{F\in \cF:x\in F\}$, $\cF_{\overline{x}}=\{F\in \cF:x\notin F\}$, $\cF_{x,y}=\{F\in \cF:x,y\in F\}$ and $\cF_{x,\overline{y}}=\{F\in \cF:x\in F,y\notin F\}$.

\begin{lem}\label{shiftdeg}
    For any left-compressed family $\cF\subseteq \binom{[n]}{k}$, we have
    \begin{enumerate}[label=\textit{(\alph*)}, ref=(\alph*)]
        \item 
        $\frac{n-k}{k}|\cF_n|\le |\cF_{\overline{n}}|$, \label{it:Fn1}
        \item 
        $\frac{|\cF_n|}{\binom{n-1}{k-1}}\le \frac{|\cF_{n-1,\overline{n}}|}{\binom{n-2}{k-1}}$. \label{it:Fn2}
    \end{enumerate}
\end{lem}

\begin{proof}
    Inequality~\ref{it:Fn1} follows from counting the pairs $(F,j)$ where $F\in\cF_n$ and $j\not\in F$. The number of such pairs is $(n-k)|\cF_n|$. Moreover, the set $S_{j,n}(F)$ is in $\cF$ and any set $F \in \cF_{\overline{n}}$ can be the $S_{j,n}$-image of at most $k$ sets in $\cF_n$. Hence the number of such pairs is at most $k|\cF_{\overline{n}}|$

    To see~\ref{it:Fn2}, a similar double counting argument establishes $\frac{n-k}{k-1}|\cF_{n,n-1}|\le |\cF_{n,\overline{n-1}}|$ and thus $|\cF_n|=|\cF_{n,n-1}|+|\cF_{n,\overline{n-1}}|\le \frac{n-1}{n-k}|\cF_{n,\overline{n-1}}|$. Also, applying $S_{n-1,n}$, we obtain $|\cF_{n,\overline{n-1}}|\le |\cF_{n-1,\overline{n}}|$. Therefore $$\frac{|\cF_n|}{\binom{n-1}{k-1}}\le \frac{\frac{n-1}{n-k}|\cF_{n,\overline{n-1}}|}{\frac{n-1}{n-k}\binom{n-2}{k-1}}\le\frac{|\cF_{n-1,\overline{n}}|}{\binom{n-2}{k-1}}$$
    as claimed.
\end{proof}



\section{Proofs}

\begin{proof}[Proof of Theorem \ref{stab}]
    Let $\cF\subseteq \binom{[sk]}{k}$ with $\nu(\cF)\le s-1$ and $\delta(\cF)\ge \delta\binom{sk-1}{k-1}$. Let us write $\cG=\binom{[sk]}{k}\setminus \cF$ and consider an arbitrary subfamily $\cG'\subseteq \cG$ with $|\cG'|=\frac{1}{s}\binom{sk}{k}=\binom{sk-1}{k-1}$. (By Theorem~\ref{K}, we know $|\cG|\ge \frac{1}{s}\binom{sk}{k}$, and so such a $\cG'$ exists.) 
    
    The minimum degree condition on $\cF$ implies $\Delta(\cG')\le (1-\delta)\binom{sk-1}{k-1}$. Corollary \ref{dp} yields that if $\delta$ is sufficiently small, then $\cG'$ (and thus $\cG$) contains at least $\frac{\delta(s-2)}{2Cs(s-1)}\binom{sk-1}{k-1}\binom{(s-1)k}{k}$ disjoint pairs. 
    
    Every set $G$ is contained in $M=\frac{1}{(s-1)!}\prod_{j=1}^{s-1}\binom{jk}{k}$ $s$-partitions of $[sk]$. Also, every disjoint pair $G,G'$ is contained in $M'=\frac{1}{(s-2)!}\prod_{j=1}^{s-2}\binom{jk}{k}$ $s$-partitions, and clearly one fixed $s$-partition contains $\binom{s}{2}$ disjoint pairs. Let us count the pairs $(G,\pi)$ with $\pi$ being an $s$-partition of $[sk]$ into $k$-sets, and $G\in \cG$ being one of those $k$-sets. On the one hand the number of such pairs is $|\cG|\cdot M$. On the other hand, by definition of $\cF$ and $\cG$, every $\pi$ contains at least one set of $\cG$. Also, by the above observations, the number of $\pi$'s containing at least two sets from $\cG$ is at least $\frac{\delta(s-2)}{2Cs(s-1)}\binom{sk-1}{k-1}\binom{(s-1)k}{k}\cdot \frac{M'}{\binom{s}{2}}$. Putting these together, 
    \[
    |\cG|\cdot M \ge \frac{1}{s!}\prod_{j=1}^{s}\binom{jk}{k}+\frac{\delta(s-2)}{2Cs(s-1)}\binom{sk-1}{k-1}\binom{(s-1)k}{k}\cdot \frac{M'}{\binom{s}{2}}.
    \]
    Dividing by $M$, we get $|\cG|\ge \binom{sk-1}{k-1}+\frac{(s-2)\delta}{s^2(s-1)C}\binom{sk-1}{k-1}$ and thus $$|\cF|= \binom{sk}{k}-|\cG|\le \binom{sk-1}{k}-\frac{(s-2)\delta}{s^2(s-1)C}\binom{sk-1}{k-1}=\Bigl(\frac{s-1}{s}-\frac{(s-2)\delta}{s^3(s-1)C}\Bigr)\binom{sk}{k}$$
    as claimed.
\end{proof}

\begin{proof}[Proof of Theorem \ref{drop}]
    Let $\cF\subseteq \binom{[sk+1]}{k}$ be a family of sets with $\nu(\cF)<s$, $|\cF|=f(sk+1,k,s)$ and let us write $n=sk+1$. 

    By Lemma \ref{shiftnu}, we can assume that $\cF$ is left-compressed. By Theorem~\ref{K}, we have $|\cF_{\overline{n}}|\le \frac{s-1}{s}\binom{sk}{k}$. Let $\delta_0$, $C$ be as in Theorem \ref{stab}, and let $\varepsilon^*_s=\min\bigl\{\frac{(s-2)\delta_0}{s^3(s-1)C},\frac{s-1}{s}-\delta_0\bigr\}$.
    
    Suppose first $|\cF_n|\le \bigl(\frac{s-1}{s}- \varepsilon^*_s\bigr)\binom{sk}{k-1}$, then $$|\cF|=|\cF_{\overline{n}}|+|\cF_n|\le \frac{s-1}{s}\binom{sk}{k}+\Bigl(\frac{s-1}{s}- \varepsilon^*_s\Bigr)\binom{sk}{k-1}\le \Bigl(\frac{s-1}{s}- \frac{\varepsilon^*_s}{s+1}\Bigr)\binom{sk+1}{k}.$$ 
    
    On the other hand, if $|\cF_n|\ge \bigl(\frac{s-1}{s}- \varepsilon^*_s\bigr)\binom{sk}{k-1}$, then we use the fact that left-compression implies $\delta(\cF_{\overline{n}})=|\cF_{n-1,\overline{n}}|$ and Lemma~\ref{shiftdeg}\ref{it:Fn2} to conclude that
    $$\frac{\delta(\cF_{\overline{n}})}{\binom{sk-1}{k-1}}=\frac{|\cF_{n-1,\overline{n}}|}{\binom{sk-1}{k-1}}\ge\frac{|\cF_{n}|}{\binom{sk}{k-1}}\ge \frac{s-1}{s}- \varepsilon^*_s .$$
    
    So one can apply Theorem~\ref{stab} with $\delta=\delta_0$ and since $\frac{s-1}{s}-\epsilon_s^*\ge\delta_0$, it gives $|\cF_{\overline{n}}|\le \bigl(\frac{s-1}{s}-\varepsilon^*_s\bigr)\binom{sk}{k}$. According to Lemma \ref{shiftdeg}\ref{it:Fn1}, we have $|\F_{n}|\le \frac{k}{(s-1)k+1}|\cF_{\overline{n}}|$ and thus
    \[
    |\cF|=|\cF_n|+|\cF_{\overline{n}}|\le \Bigl(\frac{s-1}{s}-\varepsilon^*_s\Bigr)\binom{sk+1}{k}.
    \]
    So $\varepsilon_s=\varepsilon^*_s/(s+1)$ works in both cases.
\end{proof}

\section{Concluding remarks}

As we mentioned in Section 2, Theorem \ref{drop} and Theorem \ref{stab} can be proved using the ``$\varepsilon$-$\delta$ removal lemma'' of Theorem \ref{fr} with the help of Frankl's result \cite{Fra1} on the maximum size of intersecting families with bounded maximum degree. Also, there is a large literature on the size of intersecting families with constraints on some other parameters: minimum degree \cite{HZ} and diversity \cite{Fra4,FK2,Ku,Ku2}. Supersaturation results on the minimum number of disjoint pairs in a family are also known \cite{DGS}. Although, exact results are usually known for values of $n$ large enough compared to $k$, weaker bounds could be turned into bounds on $f(n,k,s)$ with the methods used in our note.

The authors would like to acknowledge the support of Iowa State University, where this work was done while the second author was a visitor.

\end{document}